\documentclass[reqno,10pt,centertags]{amsart}
\usepackage{amsmath,amsthm,amscd,amssymb,latexsym,upref,enumerate}
\usepackage{bbm}
\usepackage{nicefrac}
\usepackage{hyperref}
\newcommand*{\mailto}[1]{\href{mailto:#1}{\nolinkurl{#1}}}
\newcommand{\arxiv}[1]{\href{http://arxiv.org/abs/#1}{arXiv:#1}} 
 
\newcommand{\N}{{\mathbb N}}

\newcommand{\C}{{\mathbb C}}

\newcommand{\bbC}{{\mathbb{C}}}

\newcommand{\bbN}{{\mathbb{N}}}

\newcommand{\bbR}{{\mathbb{R}}}

\newcommand{\cA}{{\mathcal A}}
\newcommand{\cB}{{\mathcal B}}

\newcommand{\cD}{{\mathcal D}}

\newcommand{\cH}{{\mathcal H}}

\newcommand{\cQ}{\mathcal Q}

\newcommand{\cS}{{\mathcal S}}

\newcommand{\cV}{{\mathcal V}}

\renewcommand{\a}{{\alpha}}
\renewcommand{\b}{{\beta}}
\newcommand{\g}{{\gamma}}

\newcommand{\e}{{\varepsilon}}

\DeclareMathOperator{\rank}{rank}

\DeclareMathOperator{\dom}{dom}

\renewcommand{\Im}{\text{\rm Im}}

\newcommand{\loc}{\text{\rm{loc}}}

\newcommand{\beq}{\begin{equation}}
\newcommand{\enq}{\end{equation}}

\newcommand{\no}{\notag}
\newcommand{\lb}{\label}
\newcommand{\f}{\frac}

\newcommand{\ol}{\overline}

\newcommand{\hatt}{\widehat}
\newcommand{\bi}{\bibitem}

\let\geq\geqslant
\let\leq\leqslant

\makeatletter
\def\theequation{\@arabic\c@equation}

\allowdisplaybreaks 
\numberwithin{equation}{section}

\newtheorem{theorem}{Theorem}[section]

\newtheorem{lemma}[theorem]{Lemma}
\newtheorem{corollary}[theorem]{Corollary}
\newtheorem{hypothesis}[theorem]{Hypothesis}

\theoremstyle{remark}
\newtheorem{remark}[theorem]{Remark}

\begin{document}

\title[Inverse Spectral Problems for Schr\"odinger-Type Operators]{Inverse Spectral Problems for  
 Schr\"odinger-Type Operators with Distributional Matrix-Valued Potentials}

\author[J.\ Eckhardt]{Jonathan Eckhardt}
\address{School of Computer Science \& Informatics\\ Cardiff University\\ Queen's Buildings \\ 
5 The Parade\\ Roath \\ Cardiff CF24 3AA\\ Wales \\ UK}
\email{\mailto{j.eckhardt@cs.cardiff.ac.uk}}

\author[F.\ Gesztesy]{Fritz Gesztesy}
\address{Department of Mathematics,
University of Missouri,
Columbia, MO 65211, USA}
\email{\mailto{gesztesyf@missouri.edu}}
\urladdr{\url{http://www.math.missouri.edu/personnel/faculty/gesztesyf.html}}

\author[R.\ Nichols]{Roger Nichols}
\address{Mathematics Department, The University of Tennessee at Chattanooga, 415 EMCS Building, Dept. 6956, 615 McCallie Ave, Chattanooga, TN 37403, USA}
\email{\mailto{Roger-Nichols@utc.edu}}
\urladdr{\url{http://www.utc.edu/faculty/roger-nichols/}} 

\author[A.\ L.\ Sakhnovich]{Alexander Sakhnovich}
\address{Faculty of Mathematics\\ University of Vienna\\
Oskar-Morgenstern-Platz 1\\ 1090 Wien\\ Austria}
\email{\mailto{oleksandr.sakhnovych@univie.ac.at}}
\urladdr{\url{http://www.mat.univie.ac.at/~sakhnov/}} 

\author[G.\ Teschl]{Gerald Teschl}
\address{Faculty of Mathematics\\ University of Vienna\\
Oskar-Morgenstern-Platz 1\\ 1090 Wien\\ Austria\\ and International
Erwin Schr\"odinger
Institute for Mathematical Phy\-sics\\ Boltzmanngasse 9\\ 1090 Wien\\ Austria}
\email{\mailto{Gerald.Teschl@univie.ac.at}}
\urladdr{\url{http://www.mat.univie.ac.at/~gerald/}} 

\thanks{Research supported by the Austrian Science Fund (FWF) under Grants No.\ Y330, P24301 and J3455.
R.N. gratefully acknowledges support from an AMS--Simons Travel Grant.}
\thanks{Differential Integral Equations {\bf 28}, 505--522 (2015)} 
\date{\today}
\keywords{Inverse problems, Schr\"odinger operators, matrix-valued potentials, supersymmetry.}
\subjclass[2010]{Primary 34A55, 34B20, 34B24; Secondary 34L05, 34L40, 47A10.}

\begin{abstract} 
The principal purpose of this note is to provide a reconstruction procedure for 
 distributional matrix-valued potential coefficients of Schr\"odinger-type operators on a 
half-line from the underlying Weyl--Titchmarsh function. 
\end{abstract}

\maketitle

\section{Introduction}  \lb{s1}

This note should be viewed as an addendum to the paper \cite{EGNT14}, treating 
 distributional matrix-valued potentials for (generalized) Schr\"odinger operators based on 
an intimate connection between such Schr\"odinger operators and a particular class of 
supersymmetric Dirac-type operators, and the paper \cite{Sa14} which develops a 
reconstruction procedure for the potential coefficient of a half-line Dirac operator from 
the underlying matrix-valued Weyl--Titchmarsh function. 
As a result, we derive a constructive approach to reconstruct  
distributional matrix-valued potential coefficients of (generalized) Schr\"odinger operators on 
a half-line from the underlying matrix-valued Weyl--Titchmarsh function. The importance of 
Weyl--Titchmarsh functions in connection with inverse problems for Schr\"odinger operators, 
especially, in connection with various uniqueness-type theorems has been 
well-documented in the literature. For instance, we mention the classical two-spectra 
uniqueness results due to Borg \cite{Bo46}, \cite{Bo52}, Levinson \cite{Le49}, Levitan \cite{Le68}, 
\cite[Ch.\ 3]{Le87}, Levitan and Gasymov \cite{LG64}, Marchenko \cite{Ma73}, \cite[Ch.\ 3]{Ma11},   
(see also \cite{Do65}, \cite{Ge07}, \cite{GKM02}, \cite{GS96}, \cite{GS00}, \cite{Ma99}, \cite{Ma05}  
and the extensive lists of references therein). The constructive approach to actually reconstruct 
the potential coefficient goes well beyond uniqueness theorems and now also becomes possible 
in connection with very singular (distributional) potentials. 

For the physical relevance of matrix-valued potentials, we refer, for instance to Chadan and Sabatier 
\cite[Sect.\ XI.3, XI.4]{CS89}, Newton and Jost \cite{NJ55}, and the literature cited therein. The 
classical reference on inverse scattering for matrix-valued potentials on a half-line is Agranovich and Marchenko \cite[Ch.\ V]{AM63} (see also \cite{WK74}).   

More precisely, the half-line Dirac-type operators in $L^2([0,\infty))^{2m}$, $m \in \bbN$, 
studied in this note are of the form
\begin{align}
& (D_+ (\alpha) U)(x) = (\cD U)(x)  \text{ for a.e.\ $x > 0$,}     
\no \\
& \, U \in \dom(D_+ (\alpha)) = \big\{V \in L^2([0, \infty))^{2m} 
\, \big| \, 
V \in AC([0, R])^{2m} \, \text{for all $R > 0$};    \lb{1.1} \\
& \hspace*{6.3cm}  \, \alpha V(0)=0; \, \cD V \in 
L^2([0, \infty))^{2m} \big\},    \no
\end{align}
where the $2m \times 2m$ matrix-valued differential expression $\cD$ is given by 
\begin{equation}
\cD = \begin{pmatrix} 0 & -I_m(d/dx) + \phi(x) \\ I_m(d/dx) + \phi(x) & 0 
\end{pmatrix},     \lb{1.2}
\end{equation}
and the boundary condition parameters $\alpha\in\bbC^{m\times 2m}$ satisfy the conditions 
\begin{equation}
\alpha\alpha^*=I_m, \quad \alpha J\alpha^*=0, \, \text{  where } \, 
J = \begin{pmatrix} 0 & -I_m \\ I_m & 0 \end{pmatrix}.      \lb{1.3}
\end{equation} 
Here the $m \times m$ matrix-valued potential coefficient $\phi$ is assumed to be locally square integrable on $[0,\infty)$, that is, $\phi\in L^2([0,R])^{m\times m}$ for all $R>0$, and to satisfy $\phi(\cdot) = \phi(\cdot)^*$ a.e.\ on $[0,\infty)$.

On the other hand, we define the following two kinds of quasi-derivatives,
\begin{align}
& u^{[1,j]} (x) = u'(x) + (-1)^{j+1} \phi(x) u(x)  \text{ for a.e.\ $x > 0$,}\quad j=1,2.    
\end{align} 
Thus, introducing the $m \times m$ matrix-valued differential expressions $\tau_j$, $j=1,2$, by 
\begin{equation}
(\tau_j u)(x)
= - \big(u^{[1,j]}\big)'(x) + (-1)^{j+1}\phi(x) u^{[1,j]} (x) 
  \text{ for a.e.\ $x > 0$,}  \quad j=1,2, 
\end{equation}
one infers that formally, $\tau_j$, $j=1,2$, are of the generalized Schr\"odinger form
\begin{align}
\tau_j = - I_m \f{d^2}{dx^2} + V_j(x), \quad 
V_j(x) = \phi(x)^2 + (-1)^{j} \phi'(x), \quad j=1,2.
\end{align}  
We emphasize that while $\phi^2 \in L^1_{\loc} ([0,\infty))^{m \times m}$ represents a standard 
matrix-valued potential coefficient, in general, $\phi'$ is now a genuine distribution (unless one assumes 
in addition that $\phi \in AC_{\loc}([0,\infty))^{m \times m}$). In contrast to these half-line 
Schr\"odinger operators, the Dirac-type operators $D_+ (\alpha)$ only contain the standard potential coefficient $\phi \in L^2_{\loc}([0,\infty))^{m \times m}$. 

The differential expressions $\tau_j$ then generate the generalized half-line Schr\"odinger 
operators $H_{+,0,j}$, $j=1,2$, in $L^2([0, \infty))^m$, 
\begin{align}
& (H_{+,0,j} u)(x) = (\tau_j u)(x)  = - \big(u^{[1,j]}\big)'(x) + (-1)^{j+1}\phi(x) u^{[1,j]} (x) 
  \text{ for a.e.\ $x >0$,}    \no \\
& \, u \in \dom(H_{+,0,j}) = \big\{v \in L^2([0, \infty))^m \, \big| \, v,  v^{[1,j]} \in AC([0,R])^m 
\text{ for all $R>0$};   \no \\
& \hspace*{1.5cm}  v(0) = 0; \, \big[\big(v^{[1,j]}\big)' + (-1)^j \phi v^{[1,j]}\big] \in 
L^2([0, \infty))^m\big\}, \quad  j=1,2,    \lb{1.10}
\end{align}
the primary object studied in this note. 
 
Denoting by $M^D_+ (\, \cdot \,,\alpha)$ and $\hatt M_{+,0,j}$, $j=1,2$, the 
$m \times m$ matrix-valued Weyl--Titchmarsh functions associated to $D_+ (\alpha)$ 
and $H_{+,0,j}$, $j=1,2$, respectively, the supersymmetric approach employed in 
\cite{EGNT14} naturally leads to the fundamental identity 
\begin{equation}
\hatt M_{+,0,1} (z) = \zeta M^D_+ (\zeta,\alpha_0) 
= - z \hatt M_{+,0,2} (z)^{-1}, \quad z = \zeta^2, \; 
\zeta \in \bbC \backslash \bbR,    \lb{1.11}
\end{equation}
where $\alpha_0=(I_m \quad 0)$. 

The paper \cite{Sa14}, on the other hand, focused on the inverse spectral problem for 
half-line Dirac-type operators containing $D_+ (\alpha_0)$ as a special case, and developed 
a procedure to reconstruct the matrix-valued potential coefficient from the 
underlying $m \times m$ matrix-valued Weyl--Titchmarsh function (i.e., in our particular case 
at hand, reconstructing $\phi$ from $M^D_+ (\, \cdot \,,\alpha_0)$). The reconstruction of  $\phi$ 
from $M^D_+ (\, \cdot \,,\alpha)$ with an arbitrary $\a$ satisfying \eqref{1.3} easily follows. The 
results of \cite{Sa14} generalize earlier results obtained in \cite{Sa02} for the case of locally 
bounded potentials (see more references, historical remarks and details of the procedure in 
\cite[Ch.\ 2]{SSR13}). 

We note that generalized Schr\"odinger operators (with measure and distributional potential 
coefficients) have been studied extensively in the literature. Rather than reviewing the extensive 
literature here, we refer to \cite{ET13}, \cite{EGNT14} which contain detailed historic accounts of 
this subject. 

It remains to briefly describe the content of this paper: Section \ref{s2} recalls the basics of 
Weyl--Titchmarsh theory for half-line Dirac-type operators $D_+ (\alpha)$ and 
the generalized half-line Schr\"odinger operators $H_{+,0,j}$, $j=1,2$.
Our principal Section \ref{s4} then develops a reconstruction procedure for the $m \times m$ 
matrix-valued potential coefficient $\phi$ from the underlying $m \times m$ matrix-valued 
Weyl--Titchmarsh function $M^D_+ (\, \cdot \,,\alpha)$ and hence by \eqref{1.11} also 
for the distributional $m \times m$ matrix-valued potential coefficients 
$V_j = \phi^2 + (-1)^j \phi'$ in the generalized half-line Schr\"odinger operators $H_{+,0,j}$  
from either one of $\hatt M_{+,0,1}$ or $\hatt M_{+,0,2}$. 
For simplicity, we exclusively focus on right half-lines $[0,\infty)$ throughout this note. 
The case of left half-lines is treated in a completely analogous manner. 

Concluding, we briefly summarize some of the notation used in this paper.  
All $m\times p$ matrices $M\in\bbC^{m\times p}$ will be  considered over
the field of complex numbers $\bbC$. Moreover, $I_m$ denotes
the identity matrix in $\bbC^{m\times m}$, $M^*$
the adjoint (i.e., complex conjugate transpose), and $M^\top$ the transpose 
of the matrix $M$. 

We denote with $L^2([0,\infty))^m$ the usual space of all square integrable (with respect to the Lebesgue measure) functions on $[0,\infty)$ taking values in $\C^m$, that is,
\begin{equation}
L^2([0,\infty))^m = \bigg\{U:[0,\infty)\to\bbC^{m} \, \bigg| \, \int_0^\infty dx \,
\|U(x)\|^2_{\bbC^{m}}<\infty \bigg\},  \quad m \in \bbN.    \lb{1.8}
\end{equation}
The set of functions which are only locally square integrable on $[0,\infty)$, that is, belong to 
$L^2([0,R])^m$ for all $R>0$, will be referred to as $L^2_{\loc}([0,\infty))^m$. 
The abbreviation ``a.e.'' is employed in the contexts of ``(Lebesgue) almost every'' as well as 
``(Lebesgue) almost everywhere'' on certain sets. 

With $AC_{\loc}([0,\infty))^m$ we denote the set of all functions on $[0,\infty)$ which are locally absolutely continuous, that is, belong to $AC([0,R])^m$ for all $R>0$.
The usual Sobolev spaces will be denoted by $H^1([0,R])^m$ and their local counterpart with $H^1_{\loc}([0,\infty))^m$.
We will also encounter the space $H^{-1}_{\loc}([0,\infty))$ of distributions, which is regarded as the dual of the subspace of $H^1_0([0,\infty))$ which consists of functions with compact support in $[0,\infty)$. 
Note that this space is precisely the space of distributional derivatives of functions in $L^2_{\loc}([0,\infty))$. 

The symbol $\cB(\cH_1, \cH_2)$ denotes the Banach space of bounded operators between the Hilbert spaces $\cH_1$ and $\cH_2$, and $\cB(\cH)$ abbreviates $\cB(\cH, \cH)$.
Finally, the open complex upper half-plane is denoted by $\bbC_+ = \{z \in \bbC \,|\, \Im(z) > 0\}$.

\section{Weyl--Titchmarsh Matrices for Half-Line Dirac \\ and Schr\"odinger Operators}  \lb{s2}

In this preparatory section, we review a special case of the Weyl--Titchmarsh theory for 
half-line Dirac-type and Schr\"odinger operators discussed in detail in \cite{EGNT14}.

We start by making the following simplified assumption, when compared to \cite{EGNT14}, dictated by 
the inverse spectral approach presented in our principal Section \ref{s4}.

\begin{hypothesis} \lb{h2.1}
Suppose $\phi \in L^2_{\loc}([0,\infty))^{m \times m}$, $m\in\N$, and 
$\phi(\cdot) = \phi(\cdot)^*$ a.e.\ on $[0,\infty)$. 
\end{hypothesis}

Given Hypothesis \ref{h2.1}, we introduce the $ 2m \times 2m$ matrix-valued differential expression 
\begin{equation}
\cD = \begin{pmatrix} 0 & -I_m(d/dx) + \phi(x) \\ I_m(d/dx) + \phi(x) & 0 
\end{pmatrix}.     \lb{2.6}
\end{equation}
By \cite[Lemma\ 2.15]{CG02}, $\cD$ is in the limit point case at $\infty$. (For a subsequent and more general result we refer to \cite{LM03}, see also \cite{LM00} and \cite{LO82} for such proofs under stronger hypotheses on $\phi$). 

We emphasize that the special structure of $\cD$ in \eqref{2.6} is derived from a study of supersymmetric  Dirac-type operators in $L^2(\bbR)^{2m}$, and we refer to \cite{EGNT14} 
for a detailed treatment in this context. Furthermore, we also note that \cite{EGNT14} was inspired by 
\cite{KPST05}. 
  
In order to discuss $m \times m$ Weyl--Titchmarsh matrices corresponding to self-adjoint realizations 
of $\cD$ in $L^2([0,\infty))^{2m}$, we introduce boundary condition parameters 
$\alpha = (\alpha_1 \quad \alpha_2) \in \bbC^{m\times 2m}$ satisfying the conditions 
\begin{equation}
\alpha\alpha^*=I_m, \quad \alpha J\alpha^*=0, \, \text{  where } \, 
J = \begin{pmatrix} 0 & -I_m \\ I_m & 0 \end{pmatrix}.      \lb{2.7}
\end{equation} 
Explicitly, this reads
\begin{equation}
\alpha_1\alpha_1^* +\alpha_2\alpha_2^*=I_m, \quad \alpha_2\alpha_1^* 
-\alpha_1\alpha_2^*=0. \lb{2.8}
\end{equation}
In fact, one also has
\begin{equation}
\alpha_1^*\alpha_1 +\alpha_2^*\alpha_2=I_m, \quad \alpha_2^*\alpha_1 
-\alpha_1^*\alpha_2=0, \lb{2.9}
\end{equation}
as is clear from
\begin{equation}
\begin{pmatrix} \alpha_1 & \alpha_2\\ -\alpha_2 & \alpha_1 \end{pmatrix}
\begin{pmatrix} \alpha_1^* & -\alpha_2^*\\ \alpha_2^* & \alpha_1^*
\end{pmatrix}=I_{2m}=\begin{pmatrix} \alpha_1^* & -\alpha_2^*\\ 
\alpha_2^* & \alpha_1^* \end{pmatrix}\begin{pmatrix} \alpha_1 & \alpha_2\\
-\alpha_2 & \alpha_1 \end{pmatrix}, \lb{2.10}
\end{equation}
since any left inverse matrix is also a right inverse, and vice versa. 
Moreover, from \eqref{2.9} one obtains 
\begin{equation}
\alpha^*\alpha J + J\alpha^*\alpha =  J.    \lb{2.11}
\end{equation}
The particular choice where $\alpha$ equals 
\begin{equation}
\alpha_0=(I_m \quad 0),  
\end{equation} 
will play a fundamental role later on. 

The self-adjoint half-line Dirac operators $D_+ (\alpha)$ in 
$L^2([0, \infty))^{2m}$ associated with a
self-adjoint boundary condition at $x=0$ indexed by
$\alpha\in\bbC^{m\times 2m}$ satisfying \eqref{2.7}, are of the form 
\begin{align}
& (D_+ (\alpha) U)(x) = (\cD U)(x)  \text{ for a.e.\ $x > 0$,}     
\no \\
& \, U \in \dom(D_+ (\alpha)) = \big\{V \in L^2([0, \infty))^{2m} 
\, \big| \, 
V \in AC([0, R])^{2m} \, \text{for all $R > 0$};     \lb{2.28} \\
& \hspace*{6.3cm} \alpha V(0)=0; \, \cD V \in L^2([0, \infty))^{2m} \big\}.    \no
\end{align}

Next, we denote by $U_+ (\zeta,\,\cdot\, ,\alpha)$ the $2m \times m$ 
matrix-valued Weyl--Titchmarsh solutions of $\cD U = \zeta U$, $\zeta\in\bbC\backslash\bbR$, 
satisfying 
\begin{equation}  
U_+ (\zeta,\, \cdot \,,\alpha) \in L^2([0,\infty))^{2m \times m}, \;  
\zeta\in\bbC\backslash\bbR, 
\end{equation}
and normalized such that 
\begin{align}
U_+ (\zeta,x,\alpha)&=\begin{pmatrix}u_{+,1}(\zeta,x,\alpha) \\
u_{+,2}(\zeta,x,\alpha)  \end{pmatrix}
= \Psi(\zeta,x,\alpha)\begin{pmatrix} I_m \\
M^D_+ (\zeta,\alpha) \end{pmatrix}  \no \\
&=\begin{pmatrix}\vartheta_1(\zeta,x,\alpha)
& \varphi_1(\zeta,x,\alpha)\\
\vartheta_2(\zeta,x,\alpha)
& \varphi_2(\zeta,x,\alpha)\end{pmatrix}
\begin{pmatrix} I_m \\
M^D_+ (\zeta,\alpha) \end{pmatrix}, \quad x \geq 0.      \lb{2.13}
\end{align}
In the particular case $\alpha_0=(I_m \quad  0)$ one obtains 
\begin{equation}
U_+(\zeta,0,\alpha_0) 
= \begin{pmatrix}u_{+,1}(\zeta,0,\alpha_0) \\
u_{+,2}(\zeta,0,\alpha_0)  \end{pmatrix}
= \begin{pmatrix} I_m \\ M^D_+ (\zeta,\alpha_0) \end{pmatrix}.  
\lb{2.27}
\end{equation} 
Here $M^D_+ (\zeta,\alpha)$ represents an $m\times m$ matrix, the superscript ``$D$'' indicates the 
underlying Dirac-type operator, and the functions $\Psi(\zeta,x,\alpha)$, $\vartheta_j(\zeta,x,\alpha)$, 
and $\varphi_j(\zeta,x,\alpha)$, $j=1,2$, $\zeta \in \bbC$, are defined as follows: 
$\Psi(\zeta,\,\cdot\,,\alpha)$ satisfies $\cD \Psi = \zeta \Psi$ a.e.\ on $[0,\infty)$, normalized such that 
\begin{equation}\lb{2.14}
\Psi(\zeta,0,\alpha)=(\alpha^* \; J\alpha^*)=
\begin{pmatrix} \alpha_1^* & -\alpha_2^* \\
\alpha_2^* & \alpha_1^* \end{pmatrix}.
\end{equation}
Partitioning $\Psi(\zeta,x,\alpha)$ as follows,
\begin{equation}
\Psi(\zeta,x,\alpha) = 
\begin{pmatrix}\vartheta_1(\zeta,x,\alpha) &
\varphi_1(\zeta,x,\alpha)\\
\vartheta_2(\zeta,x,\alpha)& \varphi_2(\zeta,x,\alpha)
\end{pmatrix}, \quad \zeta \in \bbC, \; x \geq 0,      \lb{2.15}
\end{equation}
defines $\vartheta_j(\zeta,x,\alpha)$ and $\varphi_j(\zeta,x,\alpha)$, 
$j=1,2$, as $m\times m$ matrices, entire with
respect to $\zeta\in\bbC$, and normalized according to \eqref{2.14}.

The $m\times m$ matrix-valued spectral function of the Dirac-type operator $D_+ (\alpha)$ 
then generates the measure $\Omega^D_+ (\,\cdot\, ,\alpha)$ in \eqref{2.22} below. In particular, 
the matrices $M^D_+ (\zeta,\alpha)$ represent the sought after half-line Weyl--Titchmarsh matrices associated with $D_+ (\alpha)$, whose basic properties can be summarized as follows:

\begin{theorem}
[\cite{AN76}, \cite{AD56}, \cite{Ca76}, \cite{CG02}, \cite{GT00}, \cite{HS81},
\cite{HS82}, \cite{HS84}, \cite{HS86}, \cite{KS88}] \lb{t2.3} ${}$ \\ 
Suppose Hypothesis \ref{h2.1}, let $\zeta\in\bbC\backslash\bbR$,
and denote by $\alpha, \delta\in\bbC^{m\times 2m}$ matrices satisfying \eqref{2.7}.\ Then the 
following hold: \\
$(i)$ $M^D_+ (\,\cdot\,,\alpha)$ is an $m\times m$ matrix-valued Nevanlinna--Herglotz function 
of maximal rank $m$. In particular,
\begin{align}
& \Im(M^D_+ (\zeta,\alpha)) \geq 0, \quad \zeta \in \bbC_+, \\
& M^D_+ (\overline \zeta,\alpha)=M^D_+ (\zeta,\alpha)^*, 
\lb{2.17} \\
& \rank (M^D_+ (\zeta,\alpha))=m,  \lb{2.18}\\
& \lim_{\varepsilon\downarrow 0} M^D_+ (\nu+
i\varepsilon,\alpha) \text{ exists for a.e.\ $\nu\in\bbR$},\lb{2.19}\\
& M^D_+(\zeta,\alpha) = [-\alpha J \delta^* +
\alpha\delta^* M^D_+ (\zeta,\delta)] [ \alpha\delta^*
+  \alpha J \delta^*M^D_+ (\zeta,\delta)]^{-1}.   \lb{2.20} 
\end{align}
Local singularities of $M^D_+ (\,\cdot\,,\alpha)$ and 
$M^D_+ (\,\cdot\,,\alpha)^{-1}$ are necessarily real and at most of first
order in the sense that 
\begin{align}
&- \lim_{\epsilon\downarrow0}
\left(i\epsilon\,
M^D_+(\nu+i\epsilon,\alpha)\right) \geq 0, \quad 
 \lim_{\epsilon\downarrow0}
\big(i\epsilon \, M^D_+ (\nu+i\epsilon,\alpha)^{-1}\big)
\geq 0, \quad \nu \in \bbR. 
\lb{2.21}  
\end{align}
$(ii)$  $M^D_+ (\,\cdot\,,\alpha)$ admits the representation 
\begin{equation}
M^D_+ (\zeta,\alpha)=F_+(\alpha)+\int_\bbR
d\Omega_+^D (\nu,\alpha) \,
\big[(\nu - \zeta)^{-1}-\nu(1+\nu^2)^{-1}\big], \lb{2.22} 
\end{equation}
where
\begin{equation}
F_+ (\alpha)=F_+ (\alpha)^*, \quad \int_\bbR
\big\|d\Omega^D_+ (\nu,\alpha)\big\|_{\bbC^{m\times m}} \,
(1+\nu^2)^{-1}<\infty. 
\end{equation}
Moreover,
\begin{equation}
\Omega_+^D((\mu,\nu],\alpha)
=\lim_{\delta\downarrow
0}\lim_{\varepsilon\downarrow 0}\f1\pi
\int_{\mu+\delta}^{\nu+\delta} d\nu' \, \Im\big(M^D_+ (\nu'+i\varepsilon,\alpha)\big).
\end{equation}
$(iii)$ $\Im\big(M^D_+ (\,\cdot\,,\alpha)\big)$ satisfies  
\begin{align}
\Im\big(M^D_+ (\zeta,\alpha)\big) &= \Im(\zeta) \int_0^{\infty} dx \,
U_+ (\zeta,x,\alpha)^* U_+ (\zeta,x,\alpha)   \no \\
&= \Im(\zeta) \int_0^{\infty} dx \, \big[u_{+,1}(\zeta,x,\alpha)^* u_{+,1}(\zeta,x,\alpha)   \lb{2.25} \\
& \hspace*{2.65cm} +
u_{+,2}(\zeta,x,\alpha)^* u_{+,2}(\zeta,x,\alpha) \big].   \no 
\end{align}
\end{theorem}

While $\cD$ contains the locally square integrable $m \times m$ matrix-valued coefficient 
$\phi\in L^2_{\loc}([0,\infty))^{m\times m}$, the associated generalized half-line Schr\"odinger 
operators to be discussed next will exhibit distributional potentials and hence are outside 
the standard Weyl--Titchmarsh theory for Sturm--Liouville operators with locally integrable 
$m \times m$ matrix-valued potentials on $[0,\infty)$. The supersymmetric 
approach employed in \cite{EGNT14} made the transition from the usual 
$L^1_{\loc}$-potentials in Schr\"odinger operators to (matrix-valued) distributional 
$H^{-1}_{\loc}$-potentials (and more general situations) in an effortless manner. 
Here, due to our assumption that $\phi$ belongs to the space $L^2_{\loc}([0,\infty))^{m\times m}$, the corresponding
potential belongs to $H^{-1}_{\loc}([0,\infty))^{m\times m}$. 

To briefly describe the corresponding generalized half-line Schr\"odinger operators, we first introduce the following two kinds of quasi-derivatives,
\begin{align}
& u^{[1,1]} (x) = (A u)(x) = u'(x) + \phi(x) u(x)  \text{ for a.e.\ $x > 0$,}    \no \\
& u \in \dom(A) = \big\{v \in L^2([0,\infty))^m \,\big|\, v \in AC([0,R]) \, \text{for all $R>0$};   \lb{3.3} \\ 
& \hspace*{5.9cm} (v' + \phi v) \in L^2([0,\infty))^m \big\},    \no 
\end{align} 
and 
\begin{align}
& u^{[1,2]} (x) = - (A^+ u)(x) = u'(x) - \phi(x) u(x)  \text{ for a.e.\ $x > 0$,}   \no \\
& u \in \dom(A^+) = \big\{v \in L^2([0,\infty))^m \,\big|\, v \in AC([0,R]) \, \text{for all $R>0$}; 
 \lb{3.4} \\
& \hspace{6.1cm} (v' - \phi v) \in L^2([0,\infty))^m\big\}.    \no     
\end{align}
Thus, introducing the $m \times m$ matrix-valued differential expressions $\tau_j$, $j=1,2$, by 
\begin{equation}
(\tau_1 u)(x) = (A^+ A u)(x) 
= - \big(u^{[1,1]}\big)'(x) + \phi(x) u^{[1,1]} (x) 
  \text{ for a.e.\ $x > 0$,}   
\end{equation}
and
\begin{equation} 
(\tau_2 u)(x) = (A A^+ u)(x)  
= - \big(u^{[1,2]}\big)'(x) - \phi(x) u^{[1,2]} (x) 
  \text{ for a.e.\ $x > 0$,}    
\end{equation}
one infers that formally, $\tau_j$, $j=1,2$, are of the generalized Schr\"odinger form
\begin{align}
\tau_j = - I_m \f{d^2}{dx^2} + V_j(x), \quad 
V_j(x) = \phi(x)^2 + (-1)^{j} \phi'(x), \quad j=1,2.
\end{align}  
We emphasize that while $\phi^2 \in L^1_{\loc}([0,\infty))^{m \times m}$ represents a standard 
matrix-valued potential coefficient, in general, $\phi'$ is now a genuine distribution (unless one assumes 
in addition that $\phi \in AC_{\loc}([0,\infty))^{m \times m}$). In contrast to these half-line 
Schr\"odinger operators, the Dirac-type operators $D_+ (\alpha)$ only contain the standard potential coefficient $\phi \in L^2_{\loc}([0,\infty))^{m \times m}$. 

By inspection, the second-order initial value problems, 
\begin{align}
& ((\tau_j - z) f)(x) = g(x) \text{ for a.e.\ $x > 0$,}  \no\\
& f, f^{[1,j]} \in AC_{\loc}([0,\infty))^m, \; g \in L^1_{\loc} ([0,\infty))^m,    \lb{3.8} \\
& \, f(x_0) = c_0, \;\, f^{[1,j]}(x_0) = d_0,\, j=1,2,    \no 
\end{align} 
for some $x_0 \geq 0$, $c_0, d_0 \in \bbC$, 
are equivalent to the first-order initial value problems  
\begin{align}
& \begin{pmatrix} f(x) \\ f^{[1,j]}(x) \end{pmatrix}^{\prime} = 
\begin{pmatrix} (-1)^j \phi(x) & 1 \\ -z & (-1)^{j+1}\phi(x) \end{pmatrix} 
\begin{pmatrix} f(x) \\ f^{[1,j]}(x) \end{pmatrix} - 
\begin{pmatrix} 0 \\ g(x) \end{pmatrix}  \text{ for a.e.\ $x > x_0$,}  
\no \\
& \begin{pmatrix} f(x_0) \\ f^{[1,j]}(x_0) \end{pmatrix} 
= \begin{pmatrix} c_0 \\ d_0 \end{pmatrix},\quad j=1,2,  \lb{3.10}
\end{align}
respectively. Since by Hypothesis \ref{h2.1}, 
$\phi \in L^2_{\loc} ([0,\infty))^{m \times m}$ (in fact, already $\phi \in L^1_{\loc} ([0,\infty))^{m \times m}$ would be sufficient), the initial value problems in
\eqref{3.10} (and hence those in \eqref{3.8}) are uniquely solvable by \cite[Theorem\ 16.1]{Na68} 
(see also \cite[Theorem\ 10.1]{EE89} and \cite[Theorem\ 16.2]{Na68}). 

Next, suppose that for some $1\leq p \leq m$, $U = (u_1 \,\, u_2)^\top$ is a 
$2m \times p$ matrix-valued solution of $\cD U = \zeta U$, that is, 
\begin{align}
& u_j \in AC_{\loc} ([0,\infty))^{m \times p}, \; j=1,2,     \lb{3.12} \\
& u_1^{[1,1]} = A u_1 \in L^1_{\loc} ([0,\infty))^{m \times p},   \quad 
u_2^{[1,2]} = - A^+ u_2  \in L^1_{\loc} ([0,\infty))^{m \times p}.  \no 
\end{align}
Then, if $\zeta \neq 0$, the supersymmetric structure of $\cD$ in 
\eqref{2.6} actually implies that also 
\begin{align}
& u_1^{[1,1]} = A u_1 = \zeta u_2 \in AC_{\loc} ([0,\infty))^{m \times p},   \lb{3.13} \\
& u_2^{[1,2]} = - A^+ u_2 = - \zeta u_1 \in AC_{\loc}([0,\infty))^{m \times p},
\lb{3.14} 
\end{align}
and hence that $u_j$ are actually distributional $m \times p$ matrix-valued solutions of $\tau_j u = \zeta^2 u$, 
$j=1,2$, that is,
\begin{align}
\begin{split}
& u_j, u_j^{[1,j]} \in AC_{\loc} ([0,\infty))^{m \times p}, \quad 
\big(u_j^{[1,j]}\big)' \in L^1_{\loc} ([0,\infty))^{m \times p}, \\  
& \tau_j u_j = - \big(u_j^{[1,j]}\big)' + (-1)^{j+1}\phi u_j^{[1,j]} = \zeta^2 u_j,\quad j=1,2.  \lb{3.15}
\end{split}
\end{align}

Thus, applying the $L^2$-property \eqref{2.25} and 
\eqref{3.12}--\eqref{3.15} to the Weyl--Titchmarsh solutions 
$U_+ (\zeta,\,\cdot\,,\alpha)$ associated with the Dirac-type differential expression $\cD$, 
then shows that $u_{+,j}(\zeta,\,\cdot\,,\alpha)$ are 
Weyl--Titchmarsh solutions associated with $\tau_j$, $j=1,2$, replacing the complex energy 
parameter $\zeta$ with $z = \zeta^2$. Moreover, introducing the following fundamental system 
$s_j (z,\,\cdot\,), c_j (z,\,\cdot\,)$, $j=1,2$, of $m \times m$ matrix-valued solutions of $\tau_j u = z u$, 
$z \in \bbC$, $j=1,2$, normalized for arbitrary $z \in \bbC$ by 
\begin{align}
s_j(z, 0) &= 0, \quad \hspace*{2.3mm}  s_j^{[1,j]}(z, 0) = I_m,                          \lb{3.17}    \\
c_j(z, 0) &= I_m, \quad  c_j^{[1,j]}(z, 0) = 0, \qquad j=1,2,    \lb{3.18} 
\end{align}
one observes as usual that for fixed $x \in \bbR$, $s_j (\,\cdot\,, x), c_j (\,\cdot\,, x)$, $j=1,2$ 
are entire. The connection with the solutions $\varphi_j$ and $\vartheta_j$, $j=1,2$, of 
$\cD U = \zeta U$ is given by
\begin{align}
& s_1(z, x) = \zeta^{-1} \varphi_1(\zeta,x,\alpha_0), \quad c_1(z, x) 
= \vartheta_1(\zeta,x,\alpha_0), \\
& s_2(z, x) = \zeta^{-1} \vartheta_2(\zeta,x,\alpha_0), \quad c_2(z, x) 
= \varphi_2(\zeta,x,\alpha_0), \quad z=\zeta^2, \; x \geq 0.
\end{align}
In addition, introducing the Weyl--Titchmarsh solutions $\psi_{+,j}(z,\, \cdot \,)$ associated 
with $\tau_j$, $j=1,2$, via
\begin{align}
\psi_{+,1} (z,\,\cdot\,) &= u_{+,1} (\zeta,\,\cdot\, ,\alpha_0), 
\lb{3.21} \\ 
\psi_{+,2} (z,\,\cdot\,) &= u_{+,2} (\zeta,\,\cdot\,,\alpha_0)
M^D_+ (\zeta,\alpha_0)^{-1},  \quad 
z = \zeta^2, \; \zeta \in \bbC \backslash \bbR, \; j=1,2,    \lb{3.22} 
\end{align}
(the right-hand sides being independent of the choice of branch for $\zeta$) and the generalized 
Dirichlet-type $m \times m$ matrix-valued Weyl--Titchmarsh functions 
$\hatt M_{+,0,j}$ of $\tau_j$, 
\begin{align}
\hatt M_{+,0,1} (z) &= \zeta M^D_+ (\zeta,\alpha_0),  \lb{3.23} \\
\hatt M_{+,0,2} (z) &= -\zeta M^D_+ (\zeta,\alpha_0)^{-1}, \quad 
z = \zeta^2, \; \zeta \in \bbC\backslash\bbR,    \lb{3.24} 
\end{align}
one infers from \eqref{2.27} that
\begin{align}
\psi_{+,j} (z,\,\cdot\,) = c_j (z,\,\cdot\,) 
+ s_j(z,\,\cdot\,) \hatt M_{+,0,j} (z), \quad 
z \in \bbC\backslash [0,\infty), \; j=1,2.   \lb{3.25} 
\end{align}
Indeed, \eqref{3.25} follows from combining \eqref{2.27}, \eqref{3.13}, and 
\eqref{3.14} (for $p=m$), which in turn imply 
\begin{equation}
 \psi_{+,j} (z, 0) = I_m, \quad 
 \psi_{+,j}^{[1,j]} (z, 0) = \hatt M_{+,0,j} (z), \quad 
z \in \bbC\backslash [0,\infty), \; j=1,2   \lb{3.27}
\end{equation}
and the unique solvability of the initial value problems in \eqref{3.8}. We summarize this 
discussion in the following result proved in \cite{EGNT14}:

\begin{theorem} \lb{t3.1}
Assume Hypothesis \ref{h2.1} and let $\alpha_0 = (I_m \; 0)$. Denote by 
\begin{equation} 
U_+ (\zeta,\,\cdot\,,\alpha_0) = (u_{+,1}(\zeta,\,\cdot\,,\alpha_0) \;  
u_{+,2}(\zeta,\,\cdot\,,\alpha_0))^{\top}, \quad \zeta \in \bbC \backslash \bbR, 
\end{equation} 
the Weyl--Titchmarsh solution corresponding to $\cD$, and by $M^D_+ (\,\cdot\,,\alpha_0)$ the 
$m \times m$ matrix-valued half-line Weyl--Titchmarsh function corresponding to $\cD$. Then 
the $m \times m$ matrix-valued Weyl--Titchmarsh solutions 
associated with $\tau_j$, denoted by $\psi_{+,j} (z,\,\cdot\,)$, $j=1,2$, are given by \eqref{3.21} and 
\eqref{3.22}, and the $m \times m$ matrix-valued generalized Dirichlet-type 
Weyl--Titchmarsh functions $\hatt M_{+,0,j}$ of $\tau_j$, $j=1,2$, are given by 
\eqref{3.23} and \eqref{3.24}. In particular,
\begin{equation}
\hatt M_{+,0,1} (z) = \zeta M^D_+ (\zeta,\alpha_0) 
= - z \hatt M_{+,0,2} (z)^{-1}, \quad z = \zeta^2, \; 
\zeta \in \bbC \backslash \bbR.    \lb{3.28}
\end{equation}
\end{theorem}

The subscript ``$0$'' in $\hatt M_{+,0,j}$, $j=1,2$, indicates that these generalized 
Weyl--Titchmarsh matrices correspond to a Dirichlet boundary condition at the reference 
point $x=0$ in the corresponding generalized half-line Schr\"odinger operators $H_{+,0,j}$, 
$j=1,2$, in $L^2([0, \infty))^m$ defined by
\begin{align}
& (H_{+,0,j} u)(x) = (\tau_j u)(x)  = - \big(u^{[1,j]}\big)'(x) + (-1)^{j+1}\phi(x) u^{[1,j]} (x) 
  \text{ for a.e.\ $x > 0$,}    \no \\
& \, u \in \dom(H_{+,0,j}) = \big\{v \in L^2([0, \infty))^m \, \big| \, v,  v^{[1,j]} \in AC([0,R])^m \text{ for all $R>0$};   \no \\
& \hspace*{1.7cm}  v(0) = 0; \, \big[\big(v^{[1,j]}\big)' + (-1)^j \phi v^{[1,j]}\big] \in 
L^2([0, \infty))^m\big\}, \quad  j=1,2.    \lb{3.29}
\end{align}
(For more general Sturm--Liouville operators in the scalar case $m=1$ we refer to 
\cite{EGNT13} and the references therein.) 
The corresponding Green's function of $H_{+,0,j}$ is then of the familiar form
\begin{align}
G_{+,0,j} (z,x,x') &= (H_{+,0,j} - z I)^{-1} (x,x')   \no \\
& = \begin{cases} s_j(z,x) \psi_{+,j} (\ol z,x')^*, & x \leq x', \\
\psi_{+,j} (z,x) s_j(\ol z,x')^*, & x' \leq x, \end{cases}    \lb{3.31} \\
& \hspace*{-5.5mm} z \in \bbC\backslash [0,\infty), \; x,x' \in [0,\infty), \; j=1,2.      \no 
\end{align}
Of course, \eqref{3.21}--\eqref{3.28}, \eqref{3.31} extend as usual to all 
$z$ in the resolvent set of $H_{+,0,j}$, $j=1,2$.
 
We conclude this section by detailing some properties of $\hatt M_{+,0,j}$: First, we 
recall the fundamental identity
\begin{equation}\lb{2.49}
\Im\big(\hatt M_{+,0,j}(z) \big)=\Im(z)\int_0^{\infty}dx'\, \psi_{+,j}(z,x')^*\psi_{+,j}(z,x'), 
\quad z\in \bbC\backslash \bbR, \; j=1,2, 
\end{equation}
implying that $\hatt M_{+,0,j}$, $j=1,2$, are matrix-valued Nevanlinna--Herglotz functions.  
Moreover, one has the following result.

\begin{lemma}\lb{l4.4}
Assume Hypothesis \ref{h2.1} and denote by $\hatt M_{+,0,j}$, $j=1,2$, the generalized 
Dir\-ichlet-type $m\times m$ matrix-valued Weyl--Titchmarsh functions associated to $H_{+,0,j}$, $j=1,2$, as defined by \eqref{3.23} and \eqref{3.24}. Then 
$\hatt M_{+,0,j}$, $j=1,2$, are $m\times m$ matrix-valued 
Nevanlinna--Herglotz functions of maximal rank $m$.  In particular $($for $j=1,2$$)$,
\begin{align}
&\Im\big(\hatt M_{+,0,j}(z)\big)\geq0,\quad z\in \bbC_+,\lb{2.50}\\
&\hatt M_{+,0,j}(\overline{z})=\hatt M_{+,0,j}(z)^*, \lb{2.51}\\
&\rank\big(\hatt M_{+,0,j}(z)\big)=m, \lb{2.52}\\
&\lim_{\e\downarrow 0}\hatt M_{+,0,j}(\lambda+i\e) \, 
\text{ exists for a.e.\ $\lambda \in \bbR$.}     \lb{2.53}
\end{align}
\end{lemma}

\section{Inverse Spectral Theory for Half-Line Dirac-Type \\ and Schr\"odinger Operators}  \lb{s4}

Several equivalent forms of self-adjoint Dirac-type systems have been considered in the literature. 
In particular, the case of self-adjoint Dirac-type systems of the form 
\begin{equation}\lb{4.1}
\frac{d}{dx} \Upsilon(\zeta,x )=i (\zeta \mathfrak{S}_3 + \mathfrak{S}_3 \cV(x)) \Upsilon(\zeta,x ) \text{ for a.e. }x > 0,
\end{equation} 
where
\begin{equation} \lb{4.2}
\mathfrak{S}_3 = \begin{pmatrix} I_{m} & 0 \\ 0 & -I_{m} \end{pmatrix}, \qquad 
\cV(x)= \begin{pmatrix} 0 &\cQ(x) \\ \cQ(x)^{*} & 0\end{pmatrix},  \quad x \geq 0, 
\end{equation} 
$\cQ$ is an $m \times m$ matrix-valued function defined a.e.\ on $[0,\infty)$, and $\zeta \in \bbC$ represents the spectral parameter, was recently studied in \cite{Sa14}. The procedure described in 
\cite{Sa14} to solve the inverse spectral problem of recovering $\cQ$ from the underlying 
matrix-valued half-line Weyl--Titchmarsh function is based on the method of operator identities \cite{SSR13, SaL76, SaL99} (see also the references therein).

In the special case when 
\begin{equation} 
\cQ(x) = - \cQ(x)^* \, \text{ for a.e.\ $x > 0$,}    \lb{4.3} 
\end{equation} 
the system \eqref{4.1} is equivalent to the half-line Dirac-type system 
\begin{equation} 
(\cD U)(\zeta,x) = \zeta U(\zeta,x), \quad \cD=J\frac{d}{dx}
+ \begin{pmatrix} 0 & \phi(x) \\ \phi(x) & 0\end{pmatrix},  \quad x > 0,   \lb{4.4} 
\end{equation}
where 
\begin{equation} 
J=\begin{pmatrix} 0 & -I_m \\ I_m & 0\end{pmatrix},
\qquad \phi(x) = \phi(x)^* \, \text{ for a.e.\ $x > 0$,}    \lb{4.5} 
 \end{equation} 
studied in the first part of Section \ref{s2}.  
 
The explicit connection between systems \eqref{4.1} and \eqref{4.4} is given by the relations
\begin{align} &       \lb{4.6}
U(\zeta,x) = W \Upsilon(\zeta,x), \quad  \phi(x) = - i \cQ(x), \quad W:=\frac{1}{\sqrt{2}}
 \begin{pmatrix}-i I_m & i I_m \\  I_m & I_m
\end{pmatrix}.
 \end{align} 
Indeed, one easily verifies that
\begin{align} &       \lb{4.7}
-W^* J\ W=i \mathfrak{S}_3, \qquad - W^* \begin{pmatrix} 0 & \phi(x) \\ \phi(x) & 0\end{pmatrix} W = \cV(x),  \quad x>0, 
 \end{align} 
where $W$ is unitary (i.e., $W^*W = W W^* = I_{2m}$). 

In order to apply the results from \cite{Sa14} to the Dirac-type system \eqref{4.4}, we need some preparations. 
First, we recall the normalized fundamental $2m \times 2m$ solution 
$\Psi(\zeta, x, \a)$ of \eqref{4.4} as introduced in \eqref{2.14}, \eqref{2.15}, with $\alpha$ 
satisfying \eqref{2.7}--\eqref{2.11}. 

The $m \times m$ matrix-valued Weyl--Titchmarsh function $ M^D_+(\cdot,\a)$, of the system 
\eqref{4.4} on $[0, \infty)$ is then introduced by the relation
\begin{equation}  
\Psi(\zeta, x, \a)\begin{pmatrix} I_m  \\  M^D_+(\zeta,\a)\end{pmatrix} 
\in L^2\big([0,\infty)\big)^{2m \times m}, \quad \zeta \in \bbC_+.     \lb{4.8} 
\end{equation}
On the other hand, the fundamental solution $\hatt \Psi(\zeta,x)$ of the Dirac-type system 
\eqref{4.1} in \cite{Sa14} is normalized at $x=0$ by 
\begin{equation} 
\hatt \Psi(\zeta,0)=I_{2m}, \quad  \zeta \in \bbC,      \lb{4.9}
\end{equation}
and the corresponding Weyl--Titchmarsh matrix  $\hatt M^D$ is introduced 
in \cite[eq.\ (1.5)]{Sa14} by the relation
\begin{equation}   
\hatt \Psi(\zeta, x)\begin{pmatrix} I_m  \\  \hatt M^D(\zeta)\end{pmatrix} 
\in L^2\big([0,\infty)\big)^{2m \times m}, \quad \zeta \in \bbC_+.    \lb{4.10}
\end{equation}
In view of \eqref{4.6}, \eqref{2.14} and \eqref{4.9} one concludes that 
\begin{equation}  
\Psi(\zeta, x, \a) = W \hatt \Psi(\zeta, x) W^* \Psi(\zeta, 0, \a), 
\quad  \zeta \in \bbC, \; x \geq 0,   \lb{4.11}
\end{equation}
and one notes that according to \eqref{2.7}, the initial value $\Psi(\zeta, 0, \a)$ is unitary.
It is immediate that the unitary matrix $W^* \Psi(\zeta, 0, \a)$ is given by
\begin{align}   
W^* \Psi(\zeta, 0, \a) &= \frac{1}{\sqrt{2}} 
\begin{pmatrix} \a_2^*+i \a_1^*  & \a_1^*-i \a_2^*\\  \a_2^*-i \a_1^* & \a_1^*+ i \a_2^*\end{pmatrix}  = \begin{pmatrix}  \a_1^*-i \a_2^* & 0\\  0 & \a_1^*+ i \a_2^* \end{pmatrix} W^*,     \lb{4.12} 
\end{align}
where, according to \eqref{2.7}, one has 
\begin{equation}  
(\a_1+i \a_2)(\a_1^*-i \a_2^*)=I_m.      \lb{4.13}
\end{equation}
Taking into account \eqref{4.8} and \eqref{4.10}--\eqref{4.13}, one derives the equality
\begin{equation}    
\hatt M^D(\zeta)=(\a_1^*+ i \a_2^*)\big[M_+^D(\zeta,\a)-i I_m\big]
\big[M_+^D(\zeta, \a)+i I_m\big]^{-1}(\a_1+i \a_2), \quad \zeta \in \bbC_+,   \lb{4.14} 
\end{equation}
relating the matrix-valued Weyl--Titchmarsh functions for systems \eqref{4.1} and \eqref{4.4}. 
We note that the Weyl--Titchmarsh matrices for both systems are unique (due to the limit 
point property of $\cD$ at $\infty$) and that $\hatt M^D$ is contractive on $\bbC_+$. 

Since $\phi =-i \cQ$ (see \eqref{4.6}), using \eqref{4.14} we can now reformulate 
\cite[Theorems 1.4 and 4.4]{Sa14} for the case of the half-line Dirac systems at hand. For that purpose, we partition $\hatt \Psi(0,x)$ into
the $m \times m$ blocks $\b_1$, $\b_2$, $\g_1$, and $\g_2$:
\begin{equation}   
\hatt \Psi(0,x)=\begin{pmatrix}\b(x) \\  \g(x)\end{pmatrix}
=\begin{pmatrix}\b_1(x) & \b_2(x)\\ \g_1(x) & \g_2(x)\end{pmatrix}, \quad x \geq 0,   \lb{4.15}
\end{equation}
and recover $\phi$ from those blocks. The  properties of $\b$ and $\g$, which we give below, are essential for their recovery and follow  immediately from \eqref{4.1}, \eqref{4.2}, and
\eqref{4.9}: 
\begin{align}  
& \b(0)=\begin{pmatrix} I_m & 0\end{pmatrix}, \quad \g(0)=\begin{pmatrix} 0 & I_m\end{pmatrix}; 
\quad \b \mathfrak{S}_3 \b^*\equiv I_m, \quad \g \mathfrak{S}_3 \g^*\equiv -I_m,
 \lb{4.16} \\ 
& \b \mathfrak{S}_3 \g^*\equiv 0, \quad \b^{\prime} \mathfrak{S}_3 \b^* 
 = \g^{\prime} \mathfrak{S}_3 \g^*\equiv 0, \quad 
 \b^{\prime} \mathfrak{S}_3 \g^*= \phi.      \lb{4.17} 
\end{align}

Next, we introduce the operator of integration, $\cA_x \in \cB\big(L^2\big([0,x]\big)^m\big)$, $x>0$, by 
\begin{equation}
\big(\cA_xf\big)(y) =-i  \int_0^y f(t)dt; \quad y \in [0,x], \; f  \in L^2\big([0,x]\big)^m, 
\end{equation}
acting componentwise on $f$.

A direct application of \cite[Theorems 1.4 and 4.4]{Sa14} then implies the following inverse 
spectral result for the half-line Dirac operator $D_+(\alpha)$:

\begin{theorem} \lb{t2.1}
Assume Hypothesis \ref{h2.1} and consider the half-line Dirac-type operator $D_+(\alpha)$ in 
\eqref{2.28}, with associated Weyl--Titchmarsh matrix $M_+^D(\, \cdot \,,\a)$.  
Then $M_+^D(\, \cdot \,,\a)$ uniquely determines $\phi(\cdot)$ a.e.\ on $[0,\infty)$.

In order to explicitly recover $\phi(\cdot)$ from $M_+^D(\, \cdot \,,\a)$, one first recovers 
the $m \times m$ matrix-valued function $\Lambda(\cdot)$ via equality \eqref{4.14} and the formula 
\begin{equation} 
\Lambda(x)= (2\pi i)^{-1}e^{x \eta } \, \underset{a \to \infty}{{\mathrm{l.i.m.}}} 
\int_{-a}^a d\xi \, \frac{e^{-i x\xi}}{\xi +i \eta}\hatt M^D \left(\frac{\xi+i \eta}{2}\right), \quad x > 0,    \lb{4.18}
\end{equation} 
where $\eta>0$ is arbitrary and l.i.m.\ denotes the entrywise limit in the norm of  $L^2\big([0,\, \infty)\big)$. 
Then 
\begin{equation}
\Lambda \in H^1_{\loc}([0,\infty))^{m \times m}. 
\end{equation}
Introducing the bounded operator  
$\Pi_x \in \cB\big(\bbC^{2m}, \, L^2\big([0,x]\big)^m\big)$, $x > 0$, via 
\begin{equation}
\big(\Pi_x g\big)(\cdot)=\Lambda(\cdot)g_1+g_2, \quad x > 0, \; 
g=\begin{pmatrix} g_1  \\  g_2\end{pmatrix}, \quad g_1, \, g_2 \in \bbC^m,    \lb{4.19}
\end{equation}
the following operator identity, 
\begin{equation}
\cA_x \cS_x - \cS_x \cA_x^*=i \Pi_x \mathfrak{S}_3 \Pi_x^*, \quad x > 0,    \lb{4.20}
\end{equation}
leads to the boundedly invertible and strictly positive operator 
$\cS_x \in \cB\big(L^2\big([0,x]\big)^m\big)$, $x>0$, given by 
\begin{equation}       \label{4.22a}
 \big(\cS_{x}f\big)(y)=f(y)- \frac{1}{2} \int_0^{x} ds \int_{|y-s|}^{y+s} dt \, \Lambda^{\prime}\left(\frac{t+y-s}{2}\right)
\Lambda^{\prime}\left(\frac{t+s-y}{2}\right)^* f(s)
\end{equation} 
for $y\in[0,x]$ and every $f\in L^2([0,x])^{m}$. Moreover, 
\begin{equation}
\Pi_x^*\cS_x^{-1}\Pi_x \in AC_{\loc}([0,\infty))^{m \times m}, 
\end{equation} 
and hence one can define the Hamiltonian $H$ of the corresponding canonical system, 
\begin{equation}  
H(x) = \g(x)^* \g(x)=\frac{d}{dx} \big(\Pi_x^*\cS_x^{-1}\Pi_x\big) \, 
\text{ for a.e.\ $x > 0$.}  \lb{4.21} 
\end{equation}
Using \eqref{4.16} and \eqref{4.17}, one uniquely recovers $\g$ and $ \b$ from $H$ as 
described in Remark \ref{r2.2} below. Finally, one obtains $\phi$ via 
\begin{equation}   
\phi(x) = \b^{\prime}(x) \mathfrak{S}_3 \g(x)^* \, \text{ for a.e.\ $x>0$}. \lb{4.22} 
\end{equation}
\end{theorem}

\begin{remark} \lb{r2.2} We describe the recovery of $\b$ and $\g$ satisfying 
\eqref{4.16} and \eqref{4.17} from $H$ given by \eqref{4.21} in greater detail. First, one 
recovers $\g_2^{-1} \g_1$ via 
\begin{equation}
\g_2^{-1} \g_1= [\g_2^* \g_2]^{-1} \g_2^* \g_1 
= \left(\begin{pmatrix} 0  &  I_m\end{pmatrix}H
\begin{pmatrix} 0  \\  I_m \end{pmatrix}\right)^{-1}
\begin{pmatrix} 0  &  I_m\end{pmatrix}H
\begin{pmatrix} I_m  \\  0 \end{pmatrix}.      \lb{4.23} 
\end{equation}
Next, one recovers $\g_2$ from the differential equation and initial condition below, 
\begin{equation}   
\g_2^{\prime} = \g_2 \big(\g_2^{-1} \g_1\big)^{\prime} 
\big(\g_2^{-1} \g_1\big)^* (I_{m} - \g_2^{-1} \g_1 (\g_2^{-1} \g_1)^*\big)^{-1}, 
\quad \g_2 (0)=I_m.    \lb{4.24}
\end{equation}
Given $\g_2$ and $\g_2^{-1} \g_1$, one recovers $\g_1$ and $\g$.
Finally, one recovers $\b$ via the relations, 
\begin{align}  
& \b=\b_1\breve \b, \quad \breve \b:=\begin{pmatrix}  I_m & \g_1^* 
\big(\g_2^*\big)^{-1}\end{pmatrix},    \lb{4.25} \\ 
& \b_1^{\prime}=-\b_1\big[\breve \b^{\prime} \mathfrak{S}_3 \big(\breve \b\big)^*\big]
\big[\breve \b \mathfrak{S}_3 \big({\breve \b}\big)^*\big]^{-1}, \quad \b_1(0)=I_{m}.   \lb{4.26}
\end{align}
\end{remark}

\medskip

Next, combining \eqref{3.28} and \eqref{4.14} one also obtains (employing $\alpha_0 = (I_m \quad 0)$)
\begin{align}
\hatt M^D(\zeta) = (-1)^{j+1} \big[\hatt M_{+,0,j}(\zeta^2)-i \zeta I_m\big]  
\big[\hatt M_{+,0,j}(\zeta^2)+i \zeta I_m\big]^{-1}, 
\quad \zeta \in \bbC_+, \; j=1,2.   \lb{4.27} 
\end{align}
Thus, given $\phi$, one has actually reconstructed the distributional potential coefficients 
$V_j = \phi^2 + (-1)^j \phi'$ in the generalized half-line Schr\"odinger operators $H_{+,0,j}$, 
$j=1,2$: 

\begin{corollary} \lb{c4.3}
Assume Hypothesis \ref{h2.1} and consider the generalized half-line Schr\"o\-dinger 
operators $H_{+,0,j}$, $j=1,2$, with associated Dirichlet-type matrix-valued Weyl--Titchmarsh functions 
$\hatt M_{+,0,j}$, $j=1,2$. Then either one of $\hatt M_{+,0,1}$ and 
$\hatt M_{+,0,2}$ uniquely determines $\phi(\cdot)$ a.e.\ on $[0,\infty)$, and hence also 
$V_j = \phi^2 + (-1)^j \phi'$, $j=1,2$. 

In addition, $\phi(\cdot)$ is recovered from $\hatt M_{+,0,1}$ $($resp., 
$\hatt M_{+,0,2}$$)$ along the lines of \eqref{4.18}--\eqref{4.22} upon employing 
\eqref{4.27} on the right-hand side of \eqref{4.18}. 
\end{corollary}

For inverse spectral problems with distributional potentials in the scalar context $m=1$ 
we also refer to  \cite{EGNT13a}.
 


\end{document}